\input ./style/arxiv-general.cfg
\documentclass[aap,MSNbibl,seceqn,dvips]{arximspdf}
\makeatletter
   \@ifpackageloaded{graphicx}{}{\usepackage{graphicx}}
\makeatother

\doi{10.1214/15-AAP1108}
\volume{26}
\issue{2}
\pubyear{2016}
\firstpage{917}
\lastpage{932}
\docsubty{FLA}

\makeatletter
\newcommand{\eqref}[1]{(\ref{#1})}
\newtheorem{theorem}{Theorem}
\newtheorem{lemma}[theorem]{Lemma}

\newproclaim{remark}[theorem]{Remark}
\newproclaim{definition}[theorem]{Definition}

\newtheorem{thmm}{Theorem}

\newproclaim{defn}[thmm]{Definition}
\newproclaim{rem}[thmm]{Remark}
\def\calD{{\mathcal D}}
\newcommand{\qnmax}{q^n_{\max}}
\newcommand{\Sn}{S^n}
\newcommand{\qmax}{q_{\max}}
\newcommand{\mun}{\mu^n}
\newcommand{\pimax}{\pi_\mathtt{max}}
\newcommand{\pimin}{\pi_\mathtt{min}}
\newcommand{\E}{\mathbb{E}}
\renewcommand{\P}{\mathbb{P}}
\newcommand{\bbZ}{\mathbb{Z}}
\newcommand{\vep}{\varepsilon}

\newproclaim{remarks}[theorem]{Remarks}

\def\calZ{{\mathcal Z}}
\def\calV{{\mathcal V}}
\def\wT{\widetilde T}
\def\wh{\widehat}

\makeatother

\begin{document}
\begin{frontmatter}

\title{Cutoff for the noisy voter model}
\runtitle{Cutoff for the noisy voter model}

\begin{aug}
\author[A]{\fnms{J. Theodore}~\snm{Cox}\corref{}\ead[label=e1]{jtcox@syr.edu}\ead[label=u1,url]{http://as-cascade.syr.edu/profiles/pages/cox-theodore.html}\thanksref{T1}},
\author[B]{\fnms{Yuval}~\snm{Peres}\ead[label=e2]{peres@microsoft.com}\ead[label=u2,url]{http://research.microsoft.com/\textasciitilde peres/}}
\and
\author[C]{\fnms{Jeffrey E.}~\snm{Steif}\ead[label=e3]{steif@math.chalmers.se}\ead[label=u3,url]{http://www.math.chalmers.se/\textasciitilde steif/}\thanksref{T2}}
\runauthor{J.~T. Cox, Y. Peres and J.~E. Steif}
\thankstext{T1}{Supported by the Simons
Foundation and the National Science Foundation.}
\thankstext{T2}{Supported by the Swedish Research Council and the Knut and Alice Wallenberg
Foundation.}
\affiliation{Syracuse University, Microsoft Research, and Chalmers
University of Technology and G\"{o}teborg University}
\address[A]{J.~T. Cox\\
Department of Mathematics\\
Syracuse University\\
215 Carnegie Building\\
Syracuse, New York 13244-1150\\
USA\\
\printead{e1}\\
\printead{u1}}
\address{}
\address[B]{Y. Peres\\
Microsoft Research\\
1 Microsoft Way\\
Redmond, Washington 98052\\
USA\\
\printead{e2}\\
\printead{u2}}
\address[C]{J.~E. Steif\\
Department of Mathematics\\
Chalmers University of Technology\\
\quad and G\"{o}teborg University\\
SE-41296 Gothenburg\\
Sweden\\
\printead{e3}\\
\printead{u3}}
\end{aug}

%
\received{\smonth{9} \syear{2014}}
%
\revised{\smonth{1} \syear{2015}}

%
\begin{abstract}
Given a continuous time Markov Chain $\{q(x,y)\}$ on a
finite set $S$, the associated noisy voter model is the
continuous time Markov chain on $\{0,1\}^S$, which evolves
in the following way: (1) for each two sites $x$ and $y$ in $S$, the
state at
site $x$ changes to the value of the state at site
$y$ at rate $q(x,y)$; (2)~each site rerandomizes
its state at rate 1. We show that if there is a uniform
bound on the rates $\{q(x,y)\}$ and the corresponding
stationary distributions are \textit{almost} uniform, then the
mixing time has a sharp cutoff at time $\log|S|/2$ with a
window of order 1. Lubetzky and Sly proved cutoff with a
window of order 1 for the stochastic Ising model on
toroids; we obtain the special case of their result for
the cycle as a consequence of our result. Finally, we consider the model
on a star and demonstrate the surprising phenomenon
that the time it takes for the chain started at all ones
to become close in total variation to the chain started at all zeros is
of smaller order than the mixing time.
\end{abstract}

%
\end{frontmatter}

\section{Introduction}\label{s.introduction}
Consider a continuous time Markov chain on the finite set~$S$,
$|S|\ge2$, where the rate of going from $x$ to $y$ is $q(x,y)$.
We let
$\qmax:=\max\{\sum_{y\neq x} q(x,y)\dvtx x\in S\}$ be the maximum rate that
we leave a state.

Next, $(S,q)$ yields a continuous time Markov process on
$\{0,1\}^S$ called \emph{the noisy voter model with voting mechanism $(S,q)$}
(often abbreviated \emph{the noisy voter model}) where, independently,
(1) for each two sites $x$ and $y$,
the state at site $x$ changes to the value of
the state at site $y$ at rate $q(x,y)$, and
(2) each site rerandomizes its state at rate 1.
By \textit{rerandomizes}, we mean that the state at that site switches to
1 or 0, each
with probability $1/2$, independently of everything else.
The noisy voter model was introduced by Granovsky and Madras \cite{GM}.
Denoting an element of $\{0,1\}^S$ by
$\eta=\{\eta(x)\}_{x\in S}$, one can describe this
dynamic in the following way: independently at each $x\in S$,
%
\begin{eqnarray}
\label{eq:nvmrates} %
&&0 \to 1 \mbox{ at rate }
\frac{1}2 + \sum_{y\ne
x}q(x,y)\eta(y) \qquad\mbox{if
}\eta(x)=0,
\nonumber
\\[-8pt]
\\[-8pt]
\nonumber
&& 1 \to 0 \mbox{ at rate }\frac{1}2 + \sum_{y\ne
x}q(x,y)
\bigl(1-\eta(y)\bigr) \qquad\mbox{if }\eta(x)=1.
\end{eqnarray}

Observe that whether or not $(S,q)$ is irreducible, the
corresponding noisy voter model is clearly
irreducible and hence has a unique stationary distribution.
If there were no rerandomization, this would simply be
the ordinary voter model associated to $q$, which has, in the case
where $q$ is irreducible,
two absorbing states,
all 0's and all 1's. On the other hand, if there were no voter mechanism
[essentially meaning that $q(x,y)=0$ for all $x$ and $y$], then the model
would simply be continuous time random walk on the hypercube.

Throughout this paper, given $q$, we let $\{\eta_t\}_{t\ge0}$ denote
the corresponding noisy voter model, $\mu_\infty$ denote its stationary
distribution and
$\mu^{\eta}_{t}$ denote the law of $\eta_t$ when $\eta_0\equiv\eta$.
(The dependence of these on $q$ is implicit.) If we have a sequence of
such systems, we let
$\{\eta^n_t\}_{t\ge0}$, $\mu^n_\infty$ and $\mu^{n,\eta}_{t}$ denote
these objects
for the $n$th system.

Recall that the total variation distance between two probability measures
$m_1$ and $m_2$ on a finite set $\Omega$ is defined to be
\[
\| m_1- m_2\|_{\mathtt{TV}}:= \frac{1}{2} \sum
_{s \in\Omega
}\bigl|m_1(s)-m_2(s)\bigr|.
\]

Next, given a noisy voter model, for $\vep>0$, we let
\[
t_{\mathtt{mix}}(\vep):= \inf\Bigl\{t\ge0\dvtx \max_{\eta\in\{0,1\}^S} \bigl\|
\mu ^{\eta}_{t}-\mu_\infty\bigr\|_{\mathtt{TV}}\le\vep\Bigr
\}
\]
denote the $\vep$-mixing time.

The main theorem of the paper is the following.

\begin{thmm}\label{thmm:main1}
Assume that we have a sequence $(S^n,q^n)$ of continuous time Markov
chains with
$\lim_{n\to\infty}|S^n|=\infty$ and $\sup_n \qnmax< \infty$. Assume
further that
there is
$C$ such that for each $n$, there is a stationary distribution for
$(S^n,q^n)$ where the ratio of the largest and smallest point masses is
at most $C$.
(This holds, e.g., in any transitive situation.) Then, for each $\vep$,
%
\begin{equation}
\label{eq:MainMixing} t_{\mathtt{mix}}(\vep) = \tfrac{1}2 \log\bigl|S^n\bigr|
\bigl(1+o(1)\bigr).
\end{equation}

Moreover, we have that
%
\begin{equation}
\label{eq:main1} \lim_{\alpha\to\infty}\liminf_{n\to\infty} \bigl\|
\mu^{n,\mathbf1}_{({1}/2)\log|\Sn|-\alpha} - \mun_\infty\bigr\| _{\mathtt
{TV}} = 1,
\end{equation}
where \textbf{1} denotes the configuration of all 1's and
%
\begin{equation}
\label{eq:main2} \lim_{\alpha\to\infty}\limsup_{n\to\infty} \max
_{\eta\in\{0,1\}^{S^n}}\bigl\| \mu^{n,\eta}_{({1}/2)\log|\Sn|+\alpha} -
\mun_\infty\bigr\|_{\mathtt{TV}} = 0.
\end{equation}
\end{thmm}

\begin{remark*}
We will see that~\eqref{eq:main2} holds in fact whenever $\lim_{n\to
\infty}|S^n|=\infty$,
and therefore the upper bound~\eqref{eq:MainMixing} also holds under
this assumption.
\end{remark*}

Theorem~\ref{thmm:main1} tells us that under the given conditions,
the mixing time is of order $\frac{1}2\log|\Sn|$ and that there is a
cutoff with a window of size of order 1. (We define mixing
times and cutoff in Section~\ref{s.background} below.) These
assumptions are
necessary. Clearly if there is no bound on
$(\qnmax)$, then the mixing time can easily be made to be of
order 1. More interestingly, even if $(\qnmax)$ is bounded,
\eqref{eq:main1} is not necessarily true without some
condition on the set of
stationary distributions. An example of this is continuous time
random walk on the \emph{$n$-star}, which is the graph
that has one vertex with $n$ edges emanating from it.
(By \textit{continuous time random walk on a graph}, we mean that
the walker waits an exponential time and then chooses a neighbor at random.)
This will be explained in Section~\ref{s.wheel}. We also mention that
it is easy to see that the condition involving the set of
stationary distributions is not necessary
in order for~(\ref{eq:main1}) and \eqref{eq:main2} to hold since one
could take
$(\qnmax)$ going to 0 sufficiently quickly so that the voter mechanism never
comes into play.

We mention that it was proved by Ramadas \cite{R} that when
randomization occurs at any rate $\delta$, the mixing time
for the noisy voter model on any graph with $n$ vertices is
$O_{\delta}(\log n)$.

Theorem~\ref{thmm:main1} has an interesting consequence for
the stochastic Ising model on cycles. The Ising model on any graph $G=(V,E)$
with parameter (inverse temperature)
$\beta\ge0$ is the probability measure on
$\{-1,1\}^{V}$ which gives, up to a normalization factor, probability
$e^{\beta\sum_{\{x,y\}\in E}\sigma(x)\sigma(y)}$ to configuration
$\sigma$.
The stochastic Ising model on $G$ with
\textit{heat-bath} dynamics is the continuous time Markov chain on
$\{-1,1\}^{V}$ where each site at rate 1 erases its present state and
chooses to be in state
$-1$ or 1, according to the conditional distribution for the Ising model,
given the other states at
that time. For the case $(\bbZ/n\bbZ)^d$, Lubetzky and Sly (see \cite{LS1})
proved that for $d=1$ and all $\beta$, $d=2$ and all $\beta$ below the
critical value
and $d=3$ and all $\beta$ sufficiently small,
one has cutoff at some constant times $\log n$ with a window of order
$\log\log n$. In \cite{LS2}, Lubetzky and Sly improved and extended
these results
in a number of directions; in particular, they proved that the result
holds for
all $\beta$ below the critical value in all dimensions and that
the window above can be taken to be of order 1. While the arguments in
this second paper
are somehow easier, they are still quite involved, including that for $d=1$.

Interestingly, for the cycle $\bbZ/n\bbZ$, the stochastic Ising model
and the noisy voter model
(where one performs random walk on $\bbZ/n\bbZ$) turn out to be the
same model,
and hence the special case of Theorem~\ref{thmm:main1} for random walk
on the cycle is already known.
In this special case, the stochastic Ising model corresponds to the
dynamics where independently at each $x\in S^n$,
the rate at which $\sigma(x)$ flips to $-\sigma(x)$ is
%
\begin{equation}
\label{eq:Isingrates} \bigl[ 1+ \exp \bigl(2\beta\sigma(x)\bigl[\sigma(x-1) + \sigma(x+1)
\bigr] \bigr) \bigr]^{-1}.
\end{equation}
An easy calculation, which we will leave to the reader, shows that if
we consider the noisy
voter model on the cycle with $q(x,x+1)=q(x,x-1)=(e^{4\beta}-1)/4$ and
multiply time by
$\theta:= \frac{2}{1+e^{4\beta}}$, we obtain the above stochastic
Ising model.
While the work of Lubetzky and Sly implies Theorem~\ref{thmm:main1} for
the cycle (and also yields some
further results), the proof given here turns out to be easier.

Mossel and Schoenebeck \cite{MS} consider a similar type of voting
model where there is no
noise and study, among other things, the time it takes to become
absorbed. Here, properly
related to our model, they show an upper bound of order $n^3$ which
would be the correct order
for the cycle. We see, from the last part of Theorem~\ref{thmm:main1}, a
drastic change
when even small noise is introduced into the system since
now it takes only order $n\log n$ to reach equilibrium. On a related note,
Mossel and Tamuz \cite{MT} provide a fascinating survey of various
``opinion exchange dynamics.''

Earlier, we mentioned the $n$-star as providing a
counterexample to
\eqref{eq:main1} when there is no
condition imposed on the stationary distributions.
The noisy voter model on the $n$-star has an additional fascinating feature.

\begin{thmm}\label{thmm:wheel} Consider the noisy voter model
corresponding to
continuous time random walk with parameter 1 on the $n$-star with $n$
even:
\begin{longlist}[(ii)]
\item[(i)] Let ${\eta}_0$ denote any configuration which is 1 on exactly
half of the leaves.
If $n\ge3$ and $t=\frac{1}4(\log n -C)> 0$, then
%
\begin{equation}
\label{eq:wheel1}\bigl \| \mu^{\eta_0}_{t} -\mu_\infty
\bigr\|_{\mathtt{TV}} \ge\frac
{e^{C}}{48+e^{C}}.
\end{equation}
\item[(ii)] The time it takes for the distribution starting from all 1's
to be within distance $1/4$ in total variation norm from the stationary
distribution is $O(1)$.
\end{longlist}
\end{thmm}

This is quite surprising since one typically expects that for monotone systems,
the mixing time for the system should be governed by the time it takes
the two extremal
states to become close in total variation norm.

We end this \hyperref[s.introduction]{Introduction} with a brief description of the results
obtainable for
a natural version of a discrete time noisy voter model.
The input for such a model is a discrete time Markov chain on a finite
set $S$ and a parameter
$\gamma\in[0,1]$. Given these, the model is defined by first
choosing an $x$ in $S$ uniformly at random, and then
with probability $1-\gamma$, one selects $y$ with probability $P(x,y)$,
at which point
the state of $x$ changes to the state of $y$, while with probability
$\gamma$, the state
at vertex $x$ is rerandomized to be 0 or 1, each with probability
$1/2$. Discrete time analogues of~\eqref{eq:main1} [and~\eqref
{eq:main1FINITE} later on]
can easily be obtained with the exact same methods we use below.
The mixing times, however, will now be at time $\frac{|S|\log
|S|}{2\gamma}$ since we
are only updating 1 vertex at a time and rerandomizing with probability
$\gamma$.
Similarly, a~discrete time analogue of~\eqref{eq:wheel1}
can be obtained when, for example, $\gamma=1/2$; here the
relevant time will be $n\log n/2$. The connection with the Ising model
holds exactly when
moving to discrete time, but then one must consider the discrete time
version of the Ising model.
The paper by Chen and Saloff-Coste (see \cite{CS}) contains various
results which allow one
to transfer between a discrete time model and its continuous time
version (where updates are
done at the times of a Poisson process). In particular, Proposition~3.2(2) in this paper
allows us to obtain a discrete time analogue of~\eqref{eq:main2} (with
time scaled again by
$n/\gamma$) \emph{from} the continuous time version of this result.
Finally a discrete time analogue of Theorem~\ref{thmm:wheel}(ii) with
the $O(1)$ term being
replaced by an $O(n)$ term can be obtained; this is done by
modifying the proof of Theorem~20.3(ii) in \cite{LPW} to obtain a
discrete time version of
Lemma~\ref{lemma.reducedwheel} from the continuous time version of
this lemma.

The rest of the paper is organized as follows.
In Section~\ref{s.background}, we briefly recall some standard
definitions concerning
mixing times and cutoff as well as introduce some notation.
In Section~\ref{s.thmproof} we prove a stronger version of
Theorem~\ref
{thmm:main1}, namely
Theorem~\ref{thmm:main1again}.
The coalescing Markov chain descriptions of both the voter model
and the noisy voter model are important tools in its analysis. However,
in this paper,
we only need these tools for the proof of the last statement of
Theorem~\ref{thmm:main1} or
equivalently for Theorem~\ref{thmm:main1again}(ii)
(as well as in the first remark in Section~\ref{s.wheel}), and
therefore these
descriptions are discussed only at those points in the paper. Finally,
Theorem~\ref{thmm:wheel} is proved in Section~\ref{s.wheel}.

\section{Background}\label{s.background}

In this section, we recall some standard definitions.
Consider a continuous time irreducible Markov chain on a finite set
$\Omega$ with
transition matrices $\{P^t(x,y)\}_{t\ge0}$ and stationary distribution
$\pi$.
Letting $P^t(x,\cdot)$ denote the distribution at time $t$ starting
from $x$, we let
%
\begin{equation}\qquad
\label{eq:LPWdef} d(t) (x) :=\bigl \|P^t(x,\cdot)-\pi\bigr\|_{\mathtt{TV}},\qquad \bar
d(t) (x,y) := \bigl\|P^t(x,\cdot)-P^t(y,\cdot)
\bigr\|_{\mathtt{TV}}
\end{equation}
and
\[
d(t) := \max_{x\in\Omega} \,d(t) (x),\qquad \bar d(t) := \max
_{x,y\in
\Omega} \bar d(t) (x,y).
\]
Next for $\vep>0$, we let $t_{\mathtt{mix}}(\vep) := \inf\{t\ge0\dvtx d(t)\le\vep\}$ denote
the $\vep$-mixing time, and then by convention we take
$t_{\mathtt{mix}}:=t_{\mathtt{mix}}(1/4)$ and call this the \emph{mixing time}.

The following notions are very natural but are perhaps not standard.
For $\vep>0$, we also let
$t_{\mathtt{mix}}(\vep)(x) := \inf\{t\ge0\dvtx d(t)(x)\le\vep\}$
and $t_{\mathtt{mix}}(x):=t_{\mathtt{mix}}(1/4)(x)$.

Following Levin, Peres and Wilmer \cite{LPW}, we say that a sequence of
Markov chains exhibits \emph{cutoff}
if for all $\vep>0$, we have
\[
\lim_{n\to\infty} \frac{t^n_{\mathtt{mix}}(\vep)}{t^n_{\mathtt
{mix}}(1-\vep)}=1.
\]
We say that a sequence of Markov chains exhibits \emph{cutoff} with a
window of size $w_n$ if
$w_n=o(t^n_{\mathtt{mix}})$ and in addition
\[
\lim_{\alpha\to\infty}\liminf_{n\to\infty} \,d_n
\bigl(t^n_{\mathtt{mix}}-\alpha w_n\bigr)=1 \quad\mbox{and}\quad
\lim_{\alpha\to\infty}\limsup_{n\to\infty} \,d_n
\bigl(t^n_{\mathtt{mix}}+\alpha w_n\bigr)=0.
\]

For continuous time random walk with rate 1 on the hypercube of
dimension~$n$, it is known
(see \cite{DGM}) that $t^n_{\mathtt{mix}}\sim\frac{1}{4}\log n$ and
that there is cutoff
with a window of order 1. Theorem~\ref{thmm:main1} states that for the
noisy voter model,
under the given assumptions, we have
that $t^n_{\mathtt{mix}}\sim\frac{1}{2}\log n$ and that there is cutoff
with a window of order 1. (The difference of $\frac{1}{4}$ and $\frac{1}{2}$
here is simply due to the fact that continuous time random walk with
rate 1 on the
hypercube of dimension $n$ has each coordinate changing its state at
rate 1 rather than
rerandomizing at rate 1.) We point out that in most cases where cutoff
is proved, the chain
is reversible, while Theorem~\ref{thmm:main1} provides for us a large
class of
nonreversible chains.

\section{Proof of Theorem~\texorpdfstring{\protect\ref{thmm:main1}}{1}}\label{s.thmproof}

We state here a stronger and more detailed version of Theorem~\ref{thmm:main1}.
First, given any probability measure on a set, we let
\[
\pimax:=\max_{x\in S}\pi(x),\qquad \pimin:= \min_{x\in S}
\pi(x) \quad\mbox{and}\quad\rho(\pi) := \frac{\pimax}{\pimin}.
\]
Given $S$ and $q$ as above, we let $\calD(q)$ denote the collection of
stationary distributions and
let
\[
\rho(q) := \min_{\pi\in\calD(q)}\rho(\pi).
\]

\begin{thmm}\label{thmm:main1again}
\textup{(i)} Fix $S$ and $q$. Let \textbf{1} denote the configuration of all 1's and
$\alpha\ge1$, and assume that
$t:=\frac{1}2\log|S|-\alpha\ge1$. Then
%
\begin{equation}
\label{eq:main1FINITE} \bigl\| \mu^{\mathbf1}_t - \mu_\infty
\bigr\|_{\mathtt{TV}} \ge \frac{0.7e^{2\alpha}}{16(1+\qmax)^2\rho^2(q)+0.7e^{2\alpha}}.
\end{equation}
\textup{(ii)} Fix $S$ and $q$. Letting superscript $H$ denote random walk (sped down
by a factor of 2) on $\{0,1\}^S$ (i.e., $q\equiv0$), we have that for
all $t$
%
\begin{equation}
\label{eq:main2FINITE} \max_{\eta_1,\eta_2\in\{0,1\}^{S}}\bigl\|\mu^{\eta_1}_t-
\mu^{\eta
_2}_t\bigr\| _{\mathtt{TV}} \le \max_{\eta_1,\eta_2\in\{0,1\}^{S}}
\bigl\|\mu^{\eta_1,H}_t-\mu^{\eta
_2,H}_t\bigr\|
_{\mathtt{TV}}.
\end{equation}
\end{thmm}

Note that~\eqref{eq:main1FINITE} implies~\eqref{eq:main1} under the
assumptions given in
Theorem~\ref{thmm:main1}. Next, since
$\max_{\eta_1,\eta_2\in\{0,1\}^{S}}\|\mu^{\eta_1,H}_{({1}/2)\log
|S|+\alpha}-\mu^{\eta_2,H}_{({1}/2)\log|S|+\alpha}\|_{\mathtt{TV}}$
is (see \cite{DGM}) at most
$\frac{4}{\sqrt{\pi}} \int_0^{{e^{-\alpha}}/{\sqrt{8}}}
e^{-t^2} \,dt
+o(1)$ as $|S|\to\infty$,
we have that~\eqref{eq:main2FINITE} implies~\eqref{eq:main2} under the
assumption that
$\lim_{n\to\infty}|S^n|=\infty$.

\subsection{Proof of Theorem~\texorpdfstring{\protect\ref{thmm:main1again}}{3}\textup{(i)}}
\mbox{}
\begin{pf*}{Proof of Theorem~\ref{thmm:main1again}\normalfont{(i)}}
We will apply
Wilson's method for obtaining lower
bounds on mixing times; see \cite{W} or Section~13.2 in \cite{LPW}.
Choose $\pi\in\calD(q)$ which minimizes $\rho(\pi)$, and let
$\Phi(\eta):=2\sum_{x\in S} \eta(x)\pi(x)-1$. We claim that we
have that
%
\begin{equation}
\label{eq:eigenvalue} \E_\eta\bigl[\Phi(\eta_t)
\bigr]=e^{-t}\Phi(\eta).
\end{equation}
To see this, let $\eta^x$ denote the configuration
$\eta$ except that the coordinate at $x$ is changed to
$1-\eta(x)$, and note that $\Phi(\eta^x)-\Phi(\eta) =
2\pi(x)(1-2\eta(x))$. Then by~\eqref{eq:nvmrates},
\begin{eqnarray*}
&&\frac{d}{dt} \E_{\eta} \bigl(\Phi(\eta_t)
\bigr)\Big|_{t=0}
\\
&&\qquad = \sum_{x\in S} \biggl(\frac{1}2 + \sum
_{y\ne x}q(x,y)1\bigl\{\eta(y)\ne \eta(x)\bigr\} \biggr)
2\pi(x) \bigl(1-2\eta(x)\bigr)
\\
&&\qquad = -\Phi(\eta) + 2\sum_{x,  y\ne x}\pi(x)q(x,y)1\bigl\{
\eta(y)\ne \eta(x)\bigr\} \bigl(1-2\eta(x)\bigr).
\end{eqnarray*}
A calculation using the stationarity of $\pi$
shows that the last sum is zero. This proves $\frac{d}{dt} \E_{\eta}
(\Phi(\eta_t))|_{t=0}=-\Phi(\eta)$, and hence
\eqref{eq:eigenvalue} holds.

Next we claim that for any $t$,
%
\begin{equation}\qquad
\label{eq:Rbound} \E_\eta \bigl(\bigl|\Phi(\eta_{t})-\Phi(
\eta)\bigr|^2 \bigr)\le (2\pimax)^2\bigl[|S|(1+\qmax) t+
\bigl(|S|(1+\qmax)t\bigr)^2\bigr].
\end{equation}
This is because a jump of $\eta_t$ changes $\Phi$ by at most
$2\pimax$, while by \eqref{eq:nvmrates} the number of
jumps during the interval $[0,t]$ is stochastically dominated above by a
Poisson random variable with mean $|S|(1+\qmax)t$.

Now consider the \emph{discrete} time Markov chain obtained by
sampling $\eta_t$ at times which are integer multiples of
$1/|S|$. Then $\Phi$ is an eigenfunction for this discrete time chain
with eigenvalue
$\lambda:=e^{-1/|S|}\in(\frac{1}2,1)$ (if $|S|\ge2$). We can now apply
equation (13.9) from Section~13.2 of \cite{LPW} to this discrete time
Markov chain with
$t$ being $|S|(\frac{1}2\log|S|-\alpha)$,
$x$ being the configuration ${\mathbf1}$ (whose corresponding $\Phi$ value
is 1)
and $R$ being $8\pimax^2(1+\qmax)^2$; see~\eqref{eq:Rbound}. Using
$\pimax\le\rho(q)/|S|$
and multiplying the numerator and denominator of the obtained fraction
from (13.9) in~\cite{LPW}
by $|S|^2$ yields~(\ref{eq:main1FINITE});
recall our continuous time system at time $\frac{1}2\log|S|-\alpha$ is
the discrete time system at time
$|S|(\frac{1}2\log|S|-\alpha)$.
\end{pf*}

\subsection{Proof of Theorem~\texorpdfstring{\protect\ref{thmm:main1again}}{3}\normalfont{(ii)}}
For part Theorem~\ref{thmm:main1again}(ii), we need to recall for the
reader the
graphical representation for the noisy voter model in terms of
coalescing Markov chains. In preparation for this part of the proof, we
will also give a result of Evans et al.
\cite{EKPS} concerning channels for noisy trees.

%
\begin{figure}

\includegraphics{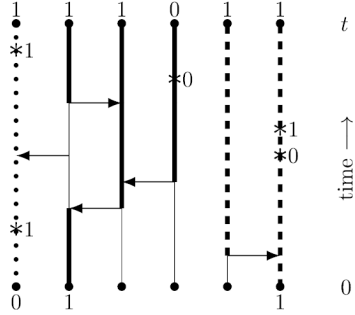}

\caption{The graphical representation and its associated
trees: arrows represent voting moves and asterisks
represent rerandomization times. In this realization,
there are three trees.}
\label{fig:harris-rep}
\end{figure}

We construct our $(S,q)$ noisy voter model using a so-called
graphical representation. Figure~\ref{fig:harris-rep} illustrates the different
elements that arise in the graphical represention. The meaning of the trees,
depicted by the dotted, solid and dashed lines will be discussed when
we get to
the proof of Theorem~\ref{thmm:main1again}(ii).
We start with the random voting times and random choices,
$T^x=\{T^x_n,n\ge1\}$ and
$W^x=\{W^x_n,n\ge1\}$, $x\in S$. The $T^x$ are independent
Poisson processes, $T^x$ has rate $q(x):=\sum_{y\ne
x}q(x,y)$ and the $W^x_n$ are independent
$S$-valued random variables, independent
of the Poisson processes, with $\P(W^x_n=y)=q(x,y)/q(x)$
for $x\ne y$. The
rerandomization times and places are given by
$R^x=\{R^x_n,n\ge1\}$ and $Z^x=\{Z^x_n,n\ge
1\}$, $x\in S$. The $R^x$ are independent rate $1$ Poisson
processes, and
the $Z^x_n$ are i.i.d.  Bernoulli random variables,
$\P(Z^x_n=1)=\P(Z^x_n=0)=1/2$.

Given $\eta_0\in\{0,1\}^S$, we define $\eta_t,t>0$ as
follows: (i) At the times $t=T^x_n$, we
draw an arrow $(x,T^x_n)\to(W^x_n,T^x_n)$ and
set $\eta_t(x) = \eta_{t-}(W^x_n)$. (ii) At
the times $t=R^x_n$, we put a $*$ at $(x,t)$ and set
$\eta_t(x)=Z^x_n$. A little thought shows that $\{\eta_t\}_{t\ge0}$ has
the dynamics specified by \eqref{eq:nvmrates}.

We construct the usual voter model dual process of coalescing
Markov chains. For
$x\in S$ and $t>0$ we construct $B^{x,t}_s,0\le
s\le t$ as follows: Set $B^{x,t}_0 =x$, and then let
$B^{x,t}_s$ trace out a path going backward in time to time
0, following the arrows for jumps. More precisely,
if $T^x\cap(0,t)=\varnothing$, put $B^{x,t}_s=x$ for $0\le s\le t$. Otherwise,
let $k=\max\{n\ge1\dvtx T^x_n<t\}$ and $u=T^x_k$, and set
\[
B^{x,t}_s = x \qquad\mbox{for }0<s<t-u \quad\mbox{and}\quad
B^{x,t}_{t-u} = W^x_k.
\]
We continue this process starting at $(B^{x,t}_{t-u},t-u)$,
thus defining $B^{x,t}_s$ for all $0\le
s\le t$. Observe that for each $x\in S$, $B^{x,t}_s$ is a
$q$-Markov chain starting at $x$. Also, these chains are
independent until they
meet, at which time they coalesce and move together thereafter.

For $t>0$, introduce $\Pi_t=\{(y,R^y_k), y\in
S,k\ge1\dvtx R^y_k\le t\}$, which contains all information up to time $t$
concerning the rerandomization \emph{times}.
For each $x\in S$, we want to look at the
time it takes the chain $B^{x,t}$ to first encounter a
rerandomization event, and also the rerandomization choice.
We do this as follows:\vspace*{1pt} If $(B^{x,t}_{s},t-s)\notin\Pi_t$
for all $0\le s\le t$, put $e(x,t)=\infty$. Otherwise,\vspace*{-1pt} let
$y,k$ satisfy $B^{x,t}_{t-R^y_k}=y$ and
$(B^{x,t}_{s},t-s)\notin\Pi_t$ for $s<t-R^y_k$, and put
$e(x,t)=t-R^y_k$ and $Z(x,t)=Z^y_k$. Given any $\eta\in
\{0,1\}^S$, the noisy voter model $\eta^\eta_t$ with
initial state $\eta^\eta_0=\eta$ can be represented as
%
\begin{equation}
\label{eq:duality} \eta^\eta_t(x) = Z(x,t)1\bigl\{e(x,t) \le t
\bigr\} + \eta\bigl(B^{x,t}_t\bigr)1\bigl\{ e(x,t)>t\bigr\},
\end{equation}
and this representation will be assumed in the rest of the
proof.

\begin{figure}[b]

\includegraphics{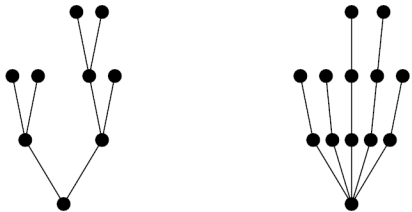}

\caption{A tree $T$ and the corresponding stringy tree
$\wT$.}\label{fig:stringy-2}
\end{figure}

In our proof of Theorem~\ref{thmm:main1again}(ii) we will use the
above graphical construction to construct certain
\emph{noisy trees} and their associated \emph{stringy trees}.
A noisy tree $T$ is a tree with flip probabilities in
$(0,\frac{1}2]$ labeling the edges. Its associated \emph{stringy tree}
$\wh{T}$
is the tree which has the same set of root--leaf paths as $T$,
but in which these paths act independently. More precisely,
for every root--leaf path in $T$, there exists an identical (in terms
of length and flip probabilities on the edges) root--leaf path in
$\wh{T}$, and in addition, all the root--leaf paths in $\wh{T}$ are
edge-disjoint. See Figure~\ref{fig:stringy-2} for an example.

Starting with $\sigma_\rho\in\{-1,+1\}$ uniform at the root
$\rho$ of $T$, we proceed upward along the tree, assigning
a value to each vertex by independently reversing the value
of the state of the parent vertex with the probability
assigned to the connecting edge (and retaining the value
otherwise). Theorem~6.1 in \cite{EKPS} relates the
conditional joint distribution (given $\sigma_\rho$) of the
resulting variables $\sigma_w$, where $w$ is a leaf of $T$
with the corresponding conditional joint distribution (given
$\sigma_{\hat{\rho}}$) for the associated stringy tree
$\wh{T}$ using \emph{channels}.
If $X$ is a random variable taking values in $\Omega_X$, and $Y$
is a random variable taking values in $\Omega_Y$,
a \emph{channel} from $X$ to $Y$ is a mapping $f\dvtx \Omega_X\times
[0,1]\rightarrow\Omega_Y$
such that if $Z$ is a uniform random variable on $[0,1]$ independent of
$X$, then
$f(X,Z)$ has distribution $Y$. See Section~15.6 in \cite{CT}.

\begin{thmm}[(Theorem~6.1 in \cite{EKPS})] \label{thmm:channel}
Given a finite noisy tree $T$ with leaves $W$ and root $\rho$,
let $\wh{T}$, with leaves
$\wh{W}$ and root $\hat{\rho}$, be the stringy tree associated with~$T$. There is a channel which, for $\xi\in\{\pm1\}$, transforms
the conditional distribution $\sigma_{\wh{W}}  |
(\sigma_{\hat{\rho}}=\xi)$
into the conditional distribution $\sigma_W   |   (\sigma_\rho=
\xi)$.
Equivalently, we say that $\wh{T}$ dominates $T$.
\end{thmm}

\begin{pf*}{Sketch of proof}
Our sketch of proof is motivated by and very similar to the
proof sketch given in \cite{P}.
We only establish a key special case of the theorem:
namely, that the tree $\Upsilon$ shown in
Figure~\ref{fig:uvtrees-3}
is dominated by the corresponding stringy tree $\wh{\Upsilon}$.
The general case is derived from it by
applying an inductive argument; see \cite{EKPS} for details.

\begin{figure}

\includegraphics{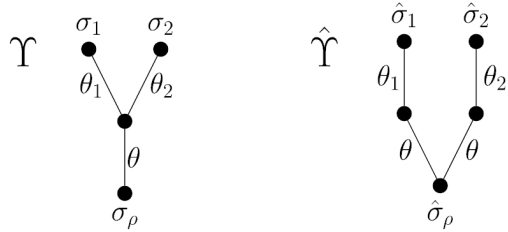}

\caption{$\Upsilon$ is dominated by $\widehat{\Upsilon}$.}
\label{fig:uvtrees-3}
\end{figure}

Let $\theta,\theta_1,\theta_2\in(0,\frac{1}2]$ be the edge flip
probabilities in Figure~\ref{fig:stringy-2}, and assume neither $\theta_1$
nor $\theta_2$ equals $\frac{1}2$ (otherwise the identity channel will
work), and w.l.o.g. assume also that $\theta_1\le\theta_2$.
Let $\sigma_\rho=\wh{\sigma}_\rho$, and let $z$ be a
$\pm1$-valued random variable, independent of
the edge flip variables, with mean
$(1-2\theta_2)/(1-2\theta_1)\in(0,1]$.
Given $0 \leq\alpha\leq1$, to be specified below,
we define the channel as follows:
%
\begin{equation}
\sigma^*_1 =\wh{\sigma}_1 \quad\mbox{and}\quad
\sigma^*_2 =\cases{ %
\wh{\sigma}_2,
&\quad $\mbox{with probability } \alpha,$
\vspace*{2pt}\cr
\wh{\sigma}_1z,& \quad$\mbox{with probability } 1-\alpha.$}
\end{equation}
It suffices to prove, for the appropriate choice of $\alpha$, that
$(\sigma_\rho,\sigma_1,\sigma_2)$
and $(\wh{\sigma}_\rho,\sigma^*_1,\sigma^*_2)$ have the same
distribution, and for this it is enough to show that the
means of all corresponding products are equal.
(This is a special case of the fact that the
characters on any finite Abelian group $G$ form a basis
for the vector space of complex functions on $G$.)
By symmetry it is only the pair correlations which require work.

Let $\gamma=1-2\theta$ and
$\gamma_i=1-2\theta_i$, $i=1,2$. Clearly
$\E(\wh{\sigma}_\rho\sigma^*_1)=\E(\sigma_\rho\sigma_1)$,
$\E(\wh{\sigma}_\rho\wh{\sigma}_1) =\gamma\gamma_1$ and
$\E(\wh{\sigma}_\rho\wh{\sigma}_2) =\gamma\gamma_2$, whence
$\E(\wh{\sigma}_\rho\sigma^*_2)
=\gamma\gamma_2=\E(\sigma_\rho\sigma_2)$ for any choice of
$\alpha$. Finally, from
$\E(\wh{\sigma}_1\wh{\sigma}_2)=\gamma^2\gamma_1\gamma_2$,
it follows that
\[
\E\bigl(\sigma^*_1\sigma^*_2\bigr) = \alpha
\gamma^2\gamma_1\gamma_2 + (1-\alpha)
\frac{\gamma_2}{\gamma_1} =\gamma_1\gamma_2 \biggl[ \alpha
\gamma^2 + (1-\alpha)\frac{1}{\gamma_1^2} \biggr].
\]
Since $\gamma^2< 1$ and $1/\gamma_1^2> 1$,
we can choose $\alpha\in[0,1]$ so that
$\E(\sigma^*_1\sigma^*_2)=\gamma_1\gamma_2=\E(\sigma_1\sigma_2)$;
explicitly,
%
\begin{equation}
\label{eq:alph} \alpha= \bigl(1-\gamma_1^2\bigr)/
\bigl(1-\gamma^2\gamma_1^2\bigr).
\end{equation}
This proves that $\wh{\Upsilon}$ dominates $\Upsilon$.
\end{pf*}

\begin{pf*}{Proof of Theorem~\ref{thmm:main1again}\normalfont{(ii)}}
Fix $t>0$ throughout.
Now for $\eta\in\{0,1\}^{S}$, consider the construction
of $\eta^{\eta}_t$ given in \eqref{eq:duality}.
Letting $\calZ(t)=\{B^{x,t}_u,x\in S, u\in[0,t]\}$, we may write
\[
\mu^{\eta}_t=\int \mu^{\eta}_{t}\bigl(
\cdot |\calZ(t)\bigr) \,d\P\bigl(\calZ(t)\bigr).
\]
Therefore, to prove \eqref{eq:main2FINITE}, it suffices to
prove the stronger fact that
for any $\eta_1,\eta_2\in\{0,1\}^{S}$ and any realization $\calZ$,
%
\begin{equation}
\label{eq:main2againagain}\bigl \| \mu^{\eta_1}_t( \cdot |\calZ)-
\mu^{\eta_2}_t( \cdot |\calZ )\bigr\|_{\mathtt{TV}} \le\max
_{\eta_1,\eta_2\in\{0,1\}^{S}}\bigl\|\mu^{\eta_1,H}_t-\mu
^{\eta
_2,H}_t\bigr\|_{\mathtt{TV}}.
\end{equation}

To proceed, we will first give for any $\eta\in S$ and
realization $\calZ$, a useful alternative description of
$\mu^{\eta}_t(  \cdot  |\calZ)$. Clearly $\calZ$ yields
a finite number of disjoint trees $T_1,T_2,\ldots,T_m$ which
describe the coalescing picture. (In the realization of
Figure~\ref{fig:harris-rep}, there are three trees indicated
by the dotted, solid and dashed lines.) Each tree has its
root sitting at $S\times\{0\}$ and its leaves sitting at $S
\times\{t\}$ in the space--time diagram. Let $x_j$ be the
root of $T_j$ and $L_j$ be the set of leaves; the $L_j$'s
are disjoint, and their union is (identified with) $S$. We
also let $\calV_j$ be the set of space--time points which consists
of the root $(x_j,0)$ along with the leaves $(\ell,t)$
and branch points of $T_j$, and view $\calV_j$ as a
tree. [If at time $s$, a chain moves from $w$ to $z$
coalescing with another walker, then we consider the branch
point to be $(z,s)$ rather than $(w,s)$.] None of this
depends on the configuration $\eta$. Note that the
branching is always at most 2 and that the tree can move
from one vertical line to another; see the solid tree in
Figure~\ref{fig:harris-rep}.

Let $Y^{\eta,j}$ be the process $\{\eta^{\eta}_s\}_{s\le t}$
conditioned on $\calZ$ restricted to $\calV_j$.
(This process also depends of course on $t$ and $\calZ$, but its dependence
on $\eta$ is what we wish to emphasize.)
Next, conditioned on $\calZ$, $Y^{\eta,1},Y^{\eta,2},\ldots,Y^{\eta,m}$
are clearly independent since $Y^{\eta,j}$ depends only on
$\Pi_t\cap T_j$ and the corresponding $Z^x_n$'s.
[This implies of course that
$\eta^{\eta}_{t}(L_1),\eta^{\eta}_{t}(L_2),\ldots,\eta^{\eta}_{t}(L_m)$
are conditionally independent given $\calZ$.]
We also let $Y^{\eta}$ be the process $\{\eta^{\eta}_s\}_{s\le t}$
restricted to $\bigcup_j \calV_j$. Crucially, $Y^{\eta,j}$ has the following
alternative simpler description as a tree-indexed Markov chain,
which is easy to verify and left to the reader.

At the root $(x_j,0)$ of $\calV_j$, the value of $Y^{\eta,j}$ is
$\eta(x_j)$.
Inductively, the value of $Y^{\eta,j}$ at a particular node is taken
to be the same as the value of its parent node
(which is lower down on the time axis) with probability
$\frac{1+e^{-s}}{2}$ where $s$ is the time difference
between these two nodes, and the opposite value otherwise.
These random choices are taken independently.
The dependence of $Y^{\eta,j}$ on $\eta$ is
only through the initial state $\eta(x_j)$; otherwise, the transition
mechanism is the same.

Consider now the process $\tilde{Y}^\eta$ indexed by $S$ and defined
by the
following two properties: the random variables $\tilde Y^\eta(x),x\in S$
are independent, and for each $j$, for all $x\in L_j$,
$\tilde{Y}^\eta(x)=\eta(x_j)$ with probability $\frac{1+e^{-t}}{2}$
and the opposite value
otherwise. It is easy to see that the distribution
of $\tilde{Y}^\eta$ is simply the distribution for
continuous time random walk on the hypercube at time $t$
started from the configuration whose state at $x$ is $\eta(x_j)$
for $x\in L_j$, $j=1,\ldots,m$.

Theorem~\ref{thmm:channel} now implies
that for each $j$, there is a channel (depending on $T_j$)
\emph{not depending on $\eta(x_j)$} which transforms the random variables
$\tilde{Y}^\eta(L_j)$ to the random variables $Y^{\eta}(L_j)=Y^{\eta
,j}(L_j)$,
meaning that given the tree $T_j$, there is a function
\[
f_{j}\dvtx \{0,1\}^{L_j}\times[0,1]\rightarrow\{0,1
\}^{L_j}
\]
so that if $U$ is a uniform random variable on $[0,1]$,
independent of everything else, we have that for each value of $\eta(x_j)$,
\[
f_{j}\bigl(\tilde{Y}^\eta(L_j), U\bigr) \quad\mbox{and}\quad {Y}^\eta(L_j)
\]
are equal in distribution.

Since $\tilde{Y}^\eta(L_j)$ are independent as we vary $j$ and similarly
for $Y^\eta(L_j)$, it follows that we have a function (depending on
$\calZ$)
\[
f\dvtx \{0,1\}^{S}\times[0,1]\rightarrow\{0,1\}^{S}
\]
so that if $U$ is a uniform random variable on $[0,1]$,
independent of everything else, we have that for any $\eta$,
\[
f\bigl(\tilde{Y}^\eta(S), U\bigr) \quad\mbox{and}\quad {Y}^\eta(S)
\]
are equal in distribution.

This then easily yields that for any $\eta_1$ and $\eta_2$,
\[
\bigl\| {Y}^{\eta_1}(S) - {Y}^{\eta_2}(S)\bigr\|_{\mathtt{TV}} \le \bigl\|
\tilde{Y}^{\eta_1}(S) - \tilde{Y}^{\eta_2}(S)\bigr\|_{\mathtt{TV}}.
\]
Finally, it is clear from construction that
\[
\bigl\| \tilde{Y}^{\eta_1}(S) - \tilde{Y}^{\eta_2}(S)\bigr\|_{\mathtt{TV}}
\le \max_{\eta_1,\eta_2\in\{0,1\}^{S}}\bigl\|\mu^{\eta_1,H}_t-
\mu^{\eta
_2,H}_t\bigr\| _{\mathtt{TV}},
\]
completing the proof.
\end{pf*}

\section{The $n$-star and the proof of Theorem~\texorpdfstring{\protect\ref{thmm:wheel}}{2}}\label
{s.wheel}

In this section, we consider the noisy voter model
$\{\eta^n_t\}$ on the $n$-star. We first explain why this
gives us an example showing that conclusion
\eqref{eq:main1} of Theorem~\ref{thmm:main1} is not true in
general without the assumption of a uniform bound on the
$\rho_n$'s even if $(\qnmax)$ is bounded. Consider first
continuous time random walk on the $n$-star with rate 1,
meaning that the walker waits an exponential amount of time
with parameter 1 and then moves to a uniform neighbor. If we
run a corresponding system of coalescing Markov chains
starting from each point, it is not hard to see that any
given pair coalesces in time $O(1)$, and that
the expected time until all chains coalesce is at most
$O(\log n)$. If we now multiply all the rates by a certain
large constant $c$, we will have that the expected time
until all chains coalesce is at most $\log n/32$. Then by
Markov's inequality, the probability that the chains have
not coalesced by time $\log n/4$ is at most $1/8$. Since
each site is rerandomized at rate 1, it is easy to see from
this fact and the graphical construction in
Section~\ref{s.thmproof} that this implies that the mixing
time at most $\log n/3$.

The rest of the section is devoted to the proof of Theorem~\ref{thmm:wheel}.

\begin{pf*}{Proof of Theorem~\ref{thmm:wheel}}
We begin with (i). This is similar to the proof of Theorem~\ref
{thmm:main1again}(i),
except one considers a different eigenfunction. Partition the leaves
into disjoint sets $A$ and $B$ each with $n/2$ elements. Let
\[
\Phi(\eta):=\sum_{x\in A}\eta(x)-\sum
_{x\in B}\eta(x).
\]
It is elementary to check that
%
\begin{equation}
\E_\eta\bigl[\Phi(\eta_t)\bigr]=e^{-2t}\Phi(
\eta).
\end{equation}
[Note that here the eigenvalue at time $t$ is $e^{-2t}$,
while in \eqref{eq:eigenvalue} it is $e^{-t}$.]

As in the proof of Theorem~\ref{thmm:main1again}(i), we
consider the discrete time Markov chain obtained by
sampling our process at times which are integer multiples of
$1/n$. Then $\Phi$ is an eigenfunction for this discrete time chain
with eigenvalue
$\lambda:=e^{-2/n}\in(\frac{1}2,1)$ (if $n\ge3$). We can now apply
equation (13.9) from Section~13.2 of \cite{LPW} to this discrete time
Markov chain with
$t$ being $\frac{n}{4}(\log n-C)$,
$x$ being the configuration $\eta_0$ (whose corresponding $\Phi$ value
is $n/2$)
and $R$ being $6$. After simplification
[and recalling that our continuous time system at time $\frac{1}4(\log
n-C)$ is the discrete time system at time
$\frac{n}{4}(\log n-C)$] we get~\eqref{eq:wheel1}.

For (ii), note first that, in the terminology introduced in
Section~\ref{s.background},
we want to show that $t^n_{\mathtt{mix}}({\mathbf1})=O(1)$.
We first note that by symmetry, if we only look at the state of the
center of the star and
the number of leaves which are in state 1, then this is also a Markov chain.
(It is a \emph{projection} of the original chain in the sense of
Section~2.3.1 in
\cite{LPW}.) Let $^{\rm R}\eta^n_t$ denote this ``reduced'' Markov chain
whose state space is $\{0,1\}\times\{0,1,\ldots, n\}$.
The key step in proving that $t^n_{\mathtt{mix}}({\mathbf1})=O(1)$
is to show that this reduced chain has mixing time $O(1)$, which is interesting
in itself; this is stated in Lemma~\ref{lemma.reducedwheel} below.

Assuming this lemma, one proceeds as follows.
Keeping symmetry in mind, we can generate a realization of the
configuration at time $t$
starting from all 1's by considering the reduced system at time $t$
starting from $(1,n)$,
and if the reduced system is in state $(a,k)$, we then construct a
configuration for
the full system by letting the center be in state $a$ and choosing a uniform
random subset of size $k$ from the $n$ leaves to be in state 1 and the
rest to be in state 0.
We can generate a realization from the stationary distribution for the
full system
in an analogous way by choosing $(a,k)$ from the stationary
distribution of the reduced system
and then letting the center be in state $a$ and choosing a uniform
random subset of size $k$ from the $n$ leaves to be in state 1 and the
rest to be in state 0.
Therefore, by an obvious coupling, we have that the total
variation distance between the full system at time $t$ started from
${\mathbf1}$ and the stationary distribution for the full system
is exactly the total variation distance
between the reduced system at time $t$ started from $(1,n)$
and the stationary distribution for the reduced system.
Now the proposition follows from Lemma~\ref{lemma.reducedwheel}.
\end{pf*}

\begin{lemma}\label{lemma.reducedwheel}
The mixing times for $\{^{\rm R}\eta^n_t\}$ is O(1).
\end{lemma}

\begin{pf}
Observe that the infinitesimal rates for this reduced chain are as follows:
\begin{eqnarray*}
(0,k)&\rightarrow&(1,k) \mbox{ at rate } \frac{1}{2}+\frac{k}{n},
\\
(0,k)&\rightarrow&(0,k+1) \mbox{ at rate } \frac{n-k}{2},
\\
(0,k)&\rightarrow&(0,k-1) \mbox{ at rate } \frac{3k}{2},
\\
(1,k)&\rightarrow&(0,k) \mbox{ at rate } \frac{1}{2}+\frac{n-k}{n},
\\
(1,k)&\rightarrow&(1,k+1) \mbox{ at rate } \frac{3(n-k)}{2},
\\
(1,k)&\rightarrow&(1,k-1) \mbox{ at rate } \frac{k}{2}.
\end{eqnarray*}
We denote this reduced system by $(X_t,Y_t)$ where $n$ is suppressed in
the notation.
The key fact that we will show is
that there exists $c_1>0$ so that for all $n$, for all (initial) states
$(a,\ell)$ and
for all (final) states $(b,k)$ with $k\in[0.4n,0.6n]$,
\[
\P_{(a,\ell)}\bigl((X_{10},Y_{10})=(b,k)\bigr)\ge
c_1/n.
\]
By equation (4.13) in \cite{LPW}, this implies that there
exists $c_2>0$ so that for all $n$,
for any two initial states, the total variation distance of
the corresponding processes at time 10 is at most $1-c_2$. This
easily leads to the claim of the lemma.

Since it is very easy for the center to change states, it is easy to
see that it suffices
to prove the above key fact when $a=1$ and $b=0$.

Let $U$ be the event that the center during $[0,10]$ never attempts an
update by looking
at one of its neighbors.
Letting $A_t:=U\cap\{X_s=1\ \forall s\in[0,t]\}$, one checks that
the conditional distribution of $Y_t$ given $A_t$
is the sum of two independent binomial distributions with
respective parameters $(\ell, \frac{3}4+\frac{1}4 e^{-2t})$ and
$(n-\ell, \frac{3}4-\frac{3}4 e^{-2t})$. In particular,
\[
g(t):=\E\biggl[\frac{Y_t}{n}\Big|A_t\biggr]=\frac{3}4+
\biggl(\frac{\ell}{n}-\frac{3}4\biggr) e^{-2t}.
\]
One also easily checks that for all $n$ and $\ell$,
%
\begin{equation}
\label{eq:gLipshitz} \bigl|g(t)-g(s)\bigr|\le2|t-s|.
\end{equation}
The representation of $Y_t$ as a sum of two binomials when conditioned
on $A_t$ yields
$\operatorname{Var}(\frac{Y_t}{n}|A_t)\le1/n$, and hence by Chebyshev's inequality we
have that
for all $n$, $\ell$, $t$ and $\sigma$,
%
\begin{equation}
\label{eq:gconcentration} \P_{(a,\ell)}\biggl(\biggl|\frac{Y_t}{n}-g(t)\biggr|\ge
\frac{\sigma}{\sqrt
{n}}\Big|A_t\biggr)\le 1/\sigma^2.
\end{equation}

Now, letting $B_t:=U\cap\{X_s=0\ \forall s\in[t,10]\}$, one checks that
the conditional distribution of $Y_{10}$ given $B_t\cap\{Y_t=nu\}$
is the sum of two independent binomial distributions with
respective parameters $(nu, \frac{1}4+\frac{3}4 e^{-2(10-t)})$ and
$(n(1-u), \frac{1}4-\frac{1}4 e^{-2(10-t)})$. In particular,
\[
h(u,t):=\E\biggl[\frac{Y_{10}}{n}\Big|B_t\cap\{Y_t=nu\}
\biggr]=\frac{1}4+\biggl(u-\frac{1}4\biggr) e^{-2(10-t)}.
\]
One also easily checks that for all $u,v$ and $t,s\in[0,10]$,
%
\begin{equation}
\label{eq:hLipshitz}\bigl |h(u,t)-h(v,s)\bigr|\le2\bigl(|u-v|+|t-s|\bigr).
\end{equation}

By an easy variant of the local central limit theorem, there exists
$c_3>0$ so that for all $n$,
$u$, $t\in[0,9.9]$ and the integers $v\in[nh(u,t)-10\sqrt
{n},nh(u,t)+10\sqrt{n}]$, one has that
%
\begin{equation}
\label{eq:LCT} \P\bigl[Y_{10}=v|B_t\cap
\{Y_t=nu\}\bigr]\ge\frac{c_3}{\sqrt{n}}.
\end{equation}

Next, one easily checks that $h(g(0),0)\le0.4$ and $h(g(9.9),9.9)\ge
0.6$, and hence by
our assumptions on $k$, there exists $t^{\star}\in[0,9.9]$ such that
$h(g(t^{\star}),t^{\star})=\frac{k}{n}$. [It is easily checked that
$h(g(t),t)$ is increasing
in $t$ but this is not needed to conclude the existence of $t^{\star}$.]

We now let $G$ be the intersection of the events $U$ and that during $[0,10]$,
the center flips its state exactly once and that this occurs during
$[t^{\star}-\frac{1}{n^{1/2}},t^{\star}+\frac{1}{n^{1/2}}]$. Clearly
there exists
$c_4>0$ so that for all $n$ and $t^{\star}$, we have that $\P(G)\ge
\frac{c_4}{\sqrt{n}}$.
On the event $G$, we let $T$ denote this unique flipping time of the center.

Now, by \eqref{eq:gLipshitz}, $|g(T)-g(t^{\star})|\le2/\sqrt{n}$
and hence
\[
\biggl\{\biggl|\frac{Y_T}{n}-g\bigl(t^{\star}\bigr)\biggr|\ge4/\sqrt{n}\biggr\}
\subseteq \biggl\{\biggl|\frac{Y_T}{n}-g(T)\biggr|\ge2/\sqrt{n}\biggr\}.
\]
Applying \eqref{eq:gconcentration}, this yields
\[
\P_{(a,\ell)}\biggl(\biggl|\frac{Y_T}{n}-g\bigl(t^{\star}\bigr)\biggr|\ge
\frac{4}{\sqrt
{n}}\Big|G,T\biggr)\le1/4.
\]
We therefore have
\[
\P_{(a,\ell)}(G\cap H)\ge\frac{c_4}{2\sqrt{n}},
\]
where $H:=\{|\frac{Y_T}{n}-g(t^{\star})|\le\frac{4}{\sqrt{n}}\}$.
Given this lower bound, to prove the key claim now, it would suffice
to show that for all parameters,
%
\begin{equation}
\label{eq:final} \P_{(a,\ell)}(Y_{10}=k|G\cap H)\ge
\frac{c_3}{\sqrt{n}},
\end{equation}
where $c_3$ comes from \eqref{eq:LCT}.

By \eqref{eq:hLipshitz}, $|T-t^{\star}| \le\frac{1}{\sqrt{n}}$ and
$|\frac{Y_T}{n}-g(t^{\star})|\le\frac{4}{\sqrt{n}}$ imply that
\[
\biggl|h\biggl(\frac{Y_T}{n},T\biggr)-h\bigl(g\bigl(t^{\star}
\bigr),t^{\star}\bigr)\biggr| \le\frac{10}{\sqrt{n}},
\]
and hence by the definition of $t^{\star}$, we have
$|h(\frac{Y_T}{n},T)-\frac{k}{n}| \le\frac{10}{\sqrt{n}}$.
Finally \eqref{eq:final} now follows from \eqref{eq:LCT} by
conditioning on
the exact values of $T$ and $Y_T$, completing the proof.
\end{pf}

\begin{remark*} In view of the proof of Theorem~\ref{thmm:wheel}(ii),
it also follows that for the reduced system, $t^n_{\mathtt{mix}}(\vep
)({\mathbf1})=O(\log(1/\vep))$.
\end{remark*}

\section*{Acknowledgments}
This work was initiated when the third author was
visiting Microsoft Research in Redmond, WA and was continued when
the first author was visiting Chalmers University of Technology.







\printaddresses
\end{document}